\documentclass[12pt,a4paper,twoside,final,notitlepage, leqno]{article}
\usepackage[english]{babel}
\usepackage[T1]{fontenc}  
\usepackage{graphicx}
\usepackage{amsmath,amsthm,epsfig,amsfonts,bbm}  

\def \smb {{\scriptstyle \bullet }}
\newcommand{\monitem}{ \smallskip \noindent $\bullet$ \quad  } 
\newcommand{\moneq}{\vspace*{-6pt} \begin{equation} \displaystyle } 
\newcommand{\moneqstar}{\vspace*{-6pt} \begin{equation*} \displaystyle } 
\newcommand{\monendstar}{\vspace*{-6pt} \end{equation*}   }
\newcommand{\monend}{\vspace*{-6pt} \end{equation}   }

\def\R{{\rm I}\! {\rm R}}

\def\br {\break}

\def\section*#1{}

\def\resume{\if@twocolumn
\section*{R\'esum\'e}
\else \small
\quotation{\bf \it R\'esum\'e \rule[1mm]{1.5mm}{0.2mm}\vspace{0pt}}
\fi}
\def\endresume{\if@twocolumn\else\endquotation\fi}
\def\abstract{\if@twocolumn
\noindent\section*{{\bf Abstract}}
\else \small
\quotation{\noindent \bf {Abstract.} \rule[1mm]{1.5mm}{0.2mm}\vspace{0pt}}
\fi}
\def\endabstract{\if@twocolumn\else\endquotation\fi}

\usepackage{fancyhdr}
\renewcommand{\headrulewidth}{0pt}
\begin{document}

\fancypagestyle{plain}{ \fancyfoot{} \renewcommand{\footrulewidth}{0pt}}
\fancypagestyle{plain}{ \fancyhead{} \renewcommand{\headrulewidth}{0pt}} 

\bibliographystyle{alpha}

\title {\bf \LARGE  Simulation of strong nonlinear waves with  \\~ 
 vectorial lattice Boltzmann schemes   \\~ }

\author { { ~  \large   Fran\c{c}ois Dubois$^{1,2}$ } \\ ~ \\   
  {\it \small $^1$    Department of Mathematics, University  Paris Sud,}\\ 
  {\it \small   B\^at. 425, F-91405 Orsay Cedex, France. }\\ 
  {\it \small $^2$   Conservatoire National des Arts et M\'etiers, Paris, France,  } \\ 
  {\it \small     Structural Mechanics and Coupled Systems Laboratory.  } \\   
  { \rm  \small francois.dubois@math.u-psud.fr. }  ~ \\ }

\bigskip

\date  {~ \\ ~ \\  {  \rm  17 November  2013} 
 \footnote {\rm  \small $\,\,$ Contribution submitted to  
  {\it   International Journal of Modern Physics C}, presented 
at the 22th  International Conference  on the Discrete Simulation of Fluid Dynamics,
 Yerevan, Armenia,   15-19 July 2013. Edition 02 January 2014.  } }

\maketitle
\renewcommand{\baselinestretch}{1.}

\markboth{Fran\c cois Dubois}
{Strong nonlinear waves with vectorial lattice Boltzmann schemes}

\bigskip

\bigskip 
\noindent  {\bf Abstract. } \qquad 
We show that an hyperbolic system with a mathematical entropy can be discretized
with vectorial lattice Boltzmann schemes with the methodology of 
kinetic representation of the dual entropy. 
We test this approach for the shallow water equations in one
and two space dimensions. We obtain interesting results for a shock tube, 
reflection of a shock  wave and unstationary two-dimensional propagation.
This contribution shows the ability of  
vectorial lattice Boltzmann schemes to simulate 
strong nonlinear  waves in unstationary situations. 
 
 $ $ \\  [2mm]

\bigskip \bigskip  \newpage \noindent {\bf \large   Introduction}  
 
\fancyfoot[C]{\oldstylenums{\thepage}}
\fancyhead[OC]{\sc{Strong nonlinear waves with vectorial lattice Boltzmann schemes}}

\monitem 
The computation of discrete shock waves with lattice Boltzmann approaches began with 
viscous Burgers approximations in the framework of lattice gaz automata
(see   Boghosian and  Levermore \cite{BL87} and  
Elton {\it et al.} \cite{ELR93}).  
With the lattice Boltzmann methods  described {\it e.g.}
 by Lallemand and Luo \cite{LL00}, first tentatives were proposed 
by  d'Humi\`eres \cite{DDH92}, Alexander 
 {\it et al.} \cite{ACCD92} among others. 
 A  D1Q2 entropic scheme for the one-dimensional viscous 
Burgers equation has been developed by Boghosian  {\it et al.} \cite{BLY04}.
The extension for gas dynamics equations 
and in particular shock tubes problems 
is studied  in   the  works of 
 Philippi {\it et al.} \cite{PHSS90},  
Nie, Shan and Chen \cite{NSC08},   
 Karlin and Asinari \cite{KA10},  
Chikatamarla and Karlin \cite{CK09}.

\monitem 
In this contribution, we experiment the ability of lattice Boltzmann schemes to approach
weak entropy  solutions of hyperbolic equations.  
It is well known that a first order hyperbolic equation exhibits shock waves. 
In order to enforce the uniqueness, the notion of mathematical entropy
has been proposed by Godunov \cite {Go59} and Friedrichs-Lax \cite{FL71}. 
A mathematical entropy is a strictly convex function of the conserved variables satisfying 
{\it ad hoc} differential constraints to ensure a complementary conservation law
for regular solutions (see {\it e.g.} our book with Despr\'es \cite{DD05}).
 The gradient of the entropy defines
the so-called ``entropy variables''. The Legendre-Fenchel-Moreau 
duality for convex functions allows us to define the dual of the entropy;
it is a convex function of the entropy variables. 

\monitem 
We start from the mathematical framework developed 
by Bouchut \cite {Bo03} making the link between the finite volume method and kinetic
models in the framework of the BGK 
approximation. 
The key notion is the representation of the dual entropy
with the help of convex functions associated with 
the discrete velocities of the lattice.  
If we suppose that a single distribution of particles is present, 
our previous contribution \cite{Du13} shows that the Burgers equation can be simulated.
We have shown also that the  approach can be extended to the nonlinear wave 
equation but is not compatible with 
the system of  shallow water equations.

\monitem 
In section~1, we develop vectorial lattice Boltzmann schemes
with kinetic  representation of the dual entropy. 
This framework is applied  in section~2 for the approximation 
of one-dimensional shallow water equations and in section~3
for the two-dimensional case. Stationary and unstationary 
 two-dimensional simulations are presented in section~4. 

%

\bigskip \bigskip  \noindent {\bf \large 1) \quad Dual entropy vectorial lattice Boltzmann schemes} 

\monitem 
In order to treat complex physics with particle like methods, 
a classical idea is to  multiply the number of particle distributions, as proposed 
by Khobalatte and Perthame~\cite{KP94}, 
Shan and Chen \cite {SC93},  
Bouchut~\cite{Bo99},        
 Dellar \cite {De02},  
  Wang {\it et al.} \cite {WWLL13}.  
We follow here  the  idea of 
dual  entropy decomposition with vectorial particle distributions  
proposed by Bouchut~\cite {Bo03}. 
We consider an hyperbolic system  composed 
by  $N$ conservation laws with  space  described by a point 
in $\, x \in \R^d$.  The 
unknowns are  the  conserved variables  $ \, W \in \R^N \,$  
 ({\it i.e.} $\, W^k   \in \R $). 
The nonlinear physical fluxes : 
  $ \,\, F_\alpha (W) \in \R^N \,$ 
(with $ \, 1 \leq \alpha \leq d $)  
are given regular  functions. The system is of first order: 
%
%
\moneq  \label{1-hyp-N}
\partial_t  W^k + \sum_{\alpha=1}^d  \partial_\alpha
 F_\alpha^{k} (W) \,=\, 0 \,, \quad 1 \leq k \leq N \, . 
\monend 
We suppose that a  mathematical entropy  $ \, \eta(W) \,$ 
is given with the associated
entropy fluxes $ \, \zeta_\alpha (W) \, $ for $\, 0 \leq \alpha \leq d $:
%
   \moneqstar 
    {\rm d} \zeta_\alpha (W) \equiv  {\rm d}  \eta(W) 
   \, \smb \, {\rm d}  F_\alpha (W) \,, \quad 1 \leq \alpha \leq d \, . 
   \monendstar 
The entropy variables  
$  \, \displaystyle \varphi_k  \equiv {{\partial \eta} (W) 
\over{\partial   W^k}}   \,$ are defined as the jacobian 
of the entropy: 
%
  \moneqstar 
  {\rm d} \eta(W) \, \equiv \, \sum_{k=1}^N  \varphi_k  \,  {\rm d}  W^k \, . 
  \monendstar 
The dual entropy   
 $   \, \eta^* ( \varphi )  \, $ 
 and the so-called ``dual entropy fluxes''  $  \, \zeta_\alpha^* ( \varphi ) \, $ 
satisfy 
\moneq \label{1-hyp-N-dual-entropy}
 \eta^* ( \varphi ) \, = \, \varphi \,\smb\, W - \eta(W) \,, \quad  
\zeta_\alpha^* ( \varphi ) \equiv \varphi \,\smb\, F_\alpha (W)
 - \zeta_\alpha (W) \, . 
\monend 
They can be differentiated without difficulty
(see {\it e.g.} \cite {DD05}):
%
   \moneqstar 
    {\rm d} \eta^*(\varphi)   \, \equiv \, 
   \sum_{k}   {\rm d} \varphi_k  \,   W^k \,, \qquad 
     {\rm d} \zeta_\alpha^*(\varphi)  \, \equiv \, 
   \sum_{k}   {\rm d} \varphi_k  \,  F_\alpha^{k} (W) \, . 
   \monendstar 
%


\monitem 
With Bouchut \cite{Bo03}, we 
introduce $N$ particule distributions  $\, f_j^k \, $ (for $  1 \leq k \leq N  $)
and $q$ velocities ($  0 \leq j \leq q-1 $).  
The conserved moments   $\, W^k   \, $  are simply the first discrete integral of these
distribution:
\moneq \label{1-hyp-N-vconserv}
 W^k \, = \,  \sum_{j=0}^{q-1} f_j^k   \, ,  \quad   1 \leq k \leq N  \, . 
\monend 
We suppose that the  particule distributions  $\, f_j^k \, $  are solution 
of the Boltzmann equations with discrete velocities:
%
\moneqstar 
   \partial_t f_j^k   + v_j^\alpha \partial_\alpha f_j^k   =  Q_j^k  
   \,\, ,   \quad 0 \leq j \leq q-1 \,, \quad  \quad  1 \leq k \leq N  
  \monendstar 
We suppose    $ \,  \sum_j  Q_j^{k}  =  0 $  in order to enforce the concervation laws
(\ref{1-hyp-N}). The nonequilibrium fluxes take the natural form
$    \Phi_\alpha^{k} \,  \equiv \, \sum_j  v_j^\alpha  f_j^k  \, $ 
and we have a system of $N$ conservation laws: 
%
   \moneqstar 
    \partial_t W^k     +  \sum_\alpha   \partial_\alpha  
     \Phi_\alpha ^{k} = 0    \, ,  \quad   1  \leq   k  \leq N  \, . 
   \monendstar

\monitem 
In the following, we call the ``Perthame-Bouchut  hypothesis'' 
the fact that  the dual  mathematical  entropy $ \, \eta^*(\varphi) \,$
can be  is decomposed into $q$   scalar potentials  $ \, h_j^* .\,$ 
The  potentials $ \, h_j^* \,$ are supposed to be regular convex functions 
of the entropy variables $ \varphi $ and satisfy the two identities 
\moneq \label{1-perthame-bouchut-hyp}
 \sum_{j=0}^{q-1}  \,  h_j^* \big(  \varphi \big)  
\,\equiv \,  \eta^* (\varphi) \,\, ,  \quad 
 \sum_{j=0}^{q-1}  \, v_j^\alpha   \,  h_j^* \big(  \varphi  \big)  \,\equiv \,  
\zeta_\alpha^* (\varphi)  \,\, ,  \quad  \forall \varphi   \, . 
\monend 
The equilibrium fluxes  
$ ( f^{\rm eq} )_j^k \,$ are easy to derive from the potentials  $ \, h_j^* \,$:
%
and  we have 
%
\moneqstar 
   \big( f^{\rm eq} \big)_j^k \,=\, {{ \partial h_j^*}\over{\partial \varphi_k}}  
   \,, \quad  
   \sum_{j=0}^{q-1} \big(  f^{\rm eq} \big)_j^{k} \,=\,    W^k \, , \quad 1 \leq k \leq N \,  
\monendstar 

  \monitem
We introduce  the Legendre dual of the convex potentials  $ \,  h_j^*  $: 
%
%
 \moneqstar 
   h_j(f_j¹ \,,\, f_j^2 \,,\, \dots \,,\, f_j^N ) \,\equiv \, \sup_{\varphi} \,
   \Big( \,  \big[ \sum_{k=1}^{N} \varphi_k \, f_j^k   \big]    -  h_j^* (\varphi) \, \Big) 
  \, , \quad 0 \leq j \leq q-1 \, . 
   \monendstar 
We observe that 
each function $\, h_j(\smb) \,$ is a  convex function  of  $N$   variables.
The so-called ``microscopic entropy''  
%
  \moneqstar    
     H(f) \equiv \sum_{j=0}^{q-1}   h_j(f_j^1 \,,\, f_j^2 \,,\, \dots \,,\, f_j^N )  \, . 
  \monendstar
It is a convex function in the domain where  the $ h_j$'s are convex.

 \monitem  
We can establish a ``H-theorem'' for the continuous dynamics relative to time
and space, 
in a way similar to the maximal entropy  approach developed by 
Karlin and his co-workers \cite {KGSB98}. 
 Under a BGK type hypothesis 
%
  \moneqstar  \label{1-bgk-hypo}
   Q_j^k \,\equiv \, {{1}\over{\tau}} \, \big( \, \big(
  f^{\rm eq} \big)_j ^k  - f_j^k  \, \big)   \,  
  \monendstar 
we have
%
 \moneqstar 
 \partial_t H(f)  + \sum_\alpha  
  \partial_\alpha \, \Big(   \sum_j   v_j^\alpha   \,  h_j(f_j^{\bullet})  
  \Big) \, \leq \, 0  \, . 
 \monendstar 

\monitem 
To establish this result, we derive the microscopic entropy relative to time: 

 \smallskip  \noindent 
$  \displaystyle   \smash { {{\partial H}\over{\partial t}} =   } \sum_{jk}   
{{ \partial h_j}\over{\partial f_j^k }} \, {{\partial  f_j^k}\over{\partial t}} $
$  \displaystyle  =    \smash { \sum_{jk}   {{ \partial h_j}\over{\partial f_j^k }} \,  Q_j^k  } $ 
$  \displaystyle -   \smash { \sum_{jk}  
{{ \partial h_j}\over{\partial f_j^k }} \, v_j^\alpha \, \partial_\alpha f_j^k  } $
$  \displaystyle =   \smash { \sum_{jk}  {{ \partial h_j}\over{\partial f_j^k }}   Q_j^k } $ 
$   \displaystyle -  \smash { \partial_\alpha \Big(  \sum_{j=0}^{q-1}   v_j^\alpha \,  h_j   \Big)  } $.   

 \smallskip  \noindent 
Then  

 \smallskip  \noindent 
  $ \, \displaystyle     \smash {  {{\partial H}\over{\partial t}}  } + 
 \partial_\alpha \big(  \sum_{j}   v_j^\alpha \,   h_j \big) \,=\, 
\, {{1}\over{\tau}} \,  \sum_{jk}   {{ \partial h_j}\over{\partial f^k }}  \big( f_j^k  \big) 
\,  \big[ \, \big( f^{\rm eq} \big)_j ^k  - f_j^k  \, \big] $ 
$ \,  \displaystyle  \smash {  \leq   {{1}\over{\tau}} \,  \sum_{jk} 
{{ \partial h_j}\over{\partial f_j^k }}  \big( f_j^{\rm eq} \big) \, 
 \big[   \big( f^{\rm eq} \big)_j ^k  - f_j^k   \big] \,}  $

 \smallskip  \noindent 
by convexity of the potentials $ \, h_j $. This last expression 
is equal to 
$ \,  \displaystyle    \smash {    {{1}\over{\tau}} \,   \sum_{jk} \, \varphi_k 
  \big[   \big( f^{\rm eq} \big)_j ^k  - f_j^k    \big] \, } $  
due to  Legendre duality: 
 \moneqstar
{{ \partial h_j}\over{\partial f_j^{\! k}}}  \big( f^{\rm eq} \big)
=  \varphi_k  \, . 
\monendstar    
%
In consequence, 
%
 \moneqstar
 {{\partial H}\over{\partial t}}  + 
 \partial_\alpha \big(  \sum_{j}   v_j^\alpha \,  h_j   \big)  \,\leq \, 
 \sum_{k} \, \varphi_k  \,    \sum_{j} \,
 \big[    \big( f^{\rm eq} \big)_j ^k  - f_j^k  \big]  =  0
\monendstar 
by construction of the values $ \, f^{\rm eq} \,$ at equilibrium. 
The H-theorem 
is proven.  \hfill $ \square $

\bigskip \bigskip  \noindent {\bf \large 2) \quad ``D1Q3Q2'' lattice Boltzmann 
scheme for shallow water} 
%

\monitem 
 We apply the previous ideas to the shollow water equations in one space 
dimension  
%
\moneqstar 
\partial_t  \rho +  \partial_x q \, =  \, 0 \,, \qquad 
\partial_t  q  +   \partial_x \Big( {{q^2}\over{\rho}} + {{p_0}\over{\rho_0^\gamma}}
  \,  \, \rho^\gamma \Big)   \, =  \, 0 \, .  
 \monendstar 
Velocity $\, u , \,$   pressure  $p$ 
and  sound velocity $ \, c>0  \, $ are given through the  expressions: 
%
\moneqstar 
     u \equiv {{q}\over{\rho}} \,, \quad 
     p \equiv  {{p_0}\over{\rho_0^\gamma}} \, \rho^\gamma  \,, \quad
     c^2 \equiv   {{\gamma \, p}\over{\rho}} \,=\,
     \gamma \, {{p_0}\over{\rho_0^\gamma}} \, \rho^{\gamma - 1} \, .  
\monendstar 
The entropy $ \, \eta \, $ and the   entropy flux   $ \, \zeta \,$  can be 
explicited without difficulty (see {\it e.g.} \cite{Du13}): 
%
\moneqstar  
     \eta \,=\, {1\over2} \, \rho \, u^2 + {{p}\over{\gamma - 1}} \,, \quad 
     \zeta \,=\, \eta \, u + p \, u \, . 
\monendstar   
Then the entropy variables   $\, \varphi = (  \theta  \equiv \partial_\rho \eta ,\, 
 \beta \equiv \partial_q  \eta ) $ 
can be related to usual ones: 
%
   \moneqstar   
   \theta = {{c^2}\over{\gamma - 1}} - {{u^2}\over{2}} \,, \quad \beta = u \, . 
   \monendstar
Thanks to (\ref{1-hyp-N-dual-entropy}), the 
dual entropy  $\, \eta^* \,$ and  the   dual entropy flux 
$\, \zeta^* \,$ can be explicited: 
$ \,  \eta^* = p $ and   $ \, \zeta^* =  p \,u $. 
We observe that
$ \, p \equiv K \,  ( \theta + {{\beta^2}\over{2}} \big) ^2 \,$
 with  $ \, K = p_0 / c_0^4 \, $   when $ \, \gamma = 2 $.
%

\monitem
We model this system with a kinetic approach and a D1Q3 stencil.
We have to find the 
particle components of the entropy variables, {\it id est} 
the (still unknown)  convex   functions $ \, h^*_j \,$ satisfing the 
Perthame-Bouchut hypothesis (\ref{1-perthame-bouchut-hyp}), 
that now can be written under the form:  
\moneq   \label{2-perthame-bouchut-st-venant} 
h_+^* ( \theta , \, \beta) + h_0^* ( \theta , \, \beta) +  h_-^* ( \theta , \, \beta)
\, = \, p \,, \quad 
\lambda \, h_+^* ( \theta , \, \beta) - \lambda \,  h_-^* ( \theta , \, \beta)
\, = \,  p \, u \,,
 \monend 
where  $\, \smash { \lambda \equiv {{\Delta x}\over{\Delta t}} } \,$ 
is the numerical velocity of the mesh. 
We use a  simple quadratic function as 
 in our previous contribution \cite{Du13}. We suggest  when $\, \gamma=2$:
\moneq   \label{2-potentiels-st-venant-1d} 
 h_0^* \,=\,  h_0^* ( \theta )  \,=\, {{a}\over{2}} \, K \, \theta^2 \, , 
 \monend 
with the introduction of a parameter $a$ that has to be precised 
for real computations. 
With this choice (\ref{2-potentiels-st-venant-1d}), 
the resolution of the system 
(\ref{2-perthame-bouchut-st-venant})  with unknowns $  h_\pm^* $ is easy:
%
\moneq   \label{2-potentiels-st-venant-1d-suite} 
 h_\pm^* (  \theta , \, \beta) \,=\, 
{{K}\over{2}} \,  \Big( \theta + {{\beta^2}\over{2}} \Big)^2 \, 
\Big( 1 \pm {{\beta}\over{\lambda}} \Big) - {{a \, K}\over{4}} \, \theta^2 \, . 
 \monend 
%

 \monitem 
From the previous potentials (\ref{2-potentiels-st-venant-1d}) and 
(\ref{2-potentiels-st-venant-1d-suite}), it is possible 
to derive the entire distribution at equilibrium. 
Observe first that with a vectorial lattice Boltzmann scheme, 
it is necessary to use 
two families  $f$ and $g$ of particle distributions, 
one relative to mass conservation and the other to momentum conservation. 
We have in this case 
%
\moneqstar 
   f_j^{\rm eq} \,=\, {{\partial h_j^*}\over{\partial \theta}} \,, \quad 
   g_j^{\rm eq} \,=\, {{\partial h_j^*}\over{\partial \beta}} \,.
\monendstar  
With (\ref{2-potentiels-st-venant-1d}), the function 
$ \, h_0^* $ is indepedent of $\beta$.  Then the particle 
$ \, \smash{ g_0 =  {{\partial h_0^*}\over{\partial \beta}} } \,$ 
is not necessary for the computation. 
With a very basic D1Q3 stencil, we define a 
``D1Q3Q2''   lattice Boltzmann scheme. 
The equilibrium distribution is obtained by 
differentiation of the relations (\ref{2-potentiels-st-venant-1d}) and 
(\ref{2-potentiels-st-venant-1d-suite}).  
%
%
From these equilibria, we implement the lattice Boltzmann method with the
MRT framework. The conserved moments follow the general paradigm 
introduced with the relation  (\ref{1-hyp-N-vconserv}): 
%
\moneqstar  
   \rho \,=\, f_0 + f_+ + f_-  \,,  \quad 
   q  \,=\, g_+ + g_-  \, . 
\monendstar
The nonconserved moments are chosen in a usual way: 
%
\moneqstar  
   J_\rho \,=\, \lambda \big( f_+ - f_- \big) \,, \quad 
   \varepsilon_\rho \,=\, \lambda^2  \big( f_+ + f_- - 2 \, f_0 \big) \,, \quad 
   J_q \,=\, \lambda \big( g_+ - g_- \big) \,. 
\monendstar
%
%
%
The relaxation step of the scheme is particularily simple when all the 
relaxation parameters  are equal to a constant value $ \, \tau \,$ as 
proposed in the BGK hypothesis. 
When a general MRT scheme is used, we follow the  rule \cite{LL00} of the 
moments $\, m_\ell^* \,$ after relaxation:
\moneq   \label{2-relaxation} 
m_\ell^* \,=\, m_\ell \,+\, s_\ell \, \big(  m_\ell^{\rm eq} -  m_\ell \big) \,. 
\monend 
%

\monitem
We have tested the previous ideas for a Riemann problem (a shock type tube).
We have chosen the following numerical data and parameters: 
%
   \moneqstar  
     \gamma =  2  \,,\,\, 
     {{\rho_\ell}\over{\rho_0}} = 2   \,,\,\, 
     {{\rho_r}\over{\rho_0}} = 0.5  \,,\,\, 
     q_\ell  =  q_r = 0    \,,\,\, 
     {{\lambda}\over{c_0}} = 8   \,,\,\, 
     a = 0.15    \,,\,\, 
     s_j \equiv  1.8 \, . 
  \monendstar  
The numerical results are displayed on Fig.~1. 
The rarefaction wave (on the left) and the 
shock wave (on the right) are correctly captured. 
%

\smallskip 
\centerline { 
\includegraphics[width=.55 \textwidth,  height=.35 \textwidth] {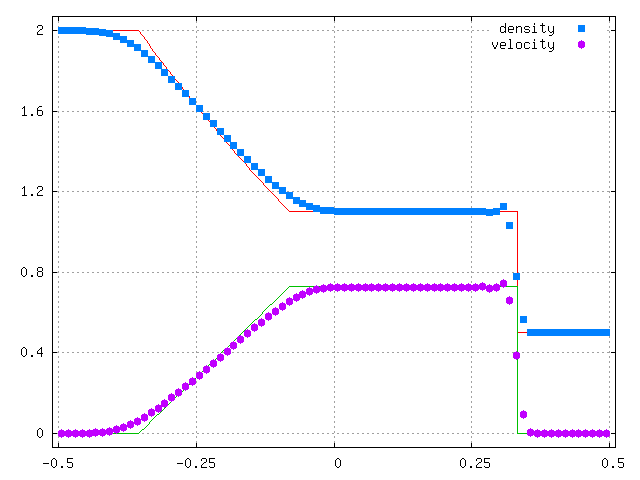}}

\noindent  {\bf Figure 1}. \quad 
 Riemann problem for shallow water equations.
 Density (blue, top) and velocity (pink, bottom) fields computed with the D1Q3Q2
lattice Boltzmann scheme with 80 mesh points and compared to the exact solution. 
\smallskip  \smallskip 

\bigskip \bigskip  \noindent {\bf \large 3) \quad ``D2Q5Q4Q4'' vectorial  lattice Boltzmann scheme} 

\monitem 
We study now the two-dimensional shallow water equations 
\moneq   \label{3-st-venant-edp} 
\left\{ \begin{array} {ll}  
\displaystyle  
\partial_t  \rho +  \partial_x \big( \rho \, u \big) 
 +  \partial_y \big( \rho \, v \big)  & \, =  \, 0 
\\ \displaystyle  
\partial_t  \big( \rho \, u \big)   
+    \partial_x \big( \rho \, u^2  + {{p_0}\over{\rho_0^2}} \, \rho^2  \big) 
+    \partial_y \big( \rho \, u \, v  \big) 
  & \, =  \, 0  
\\ \displaystyle   
\partial_t  \big( \rho \, v \big)   
+    \partial_x \big( \rho \, u \, v  \big) 
+    \partial_y \big( \rho \,  v^2  + {{p_0}\over{\rho_0^2}}  \, \rho^2   \big) 
  & \, =  \, 0 \, . 
\end{array} \right. \monend 
We have three conservation laws in two space dimensions.
We extend the previous D1Q3Q2 vectorial lattice Boltzmann scheme 
into a  D2Q5Q4Q4 scheme. 
The D2Q5 stencil is associated to the following velocities:
\moneq   \label{3-stencil-d2q5} 
 v_0 = (0,\, 0) \,,\quad 
 v_1 = (\lambda,\, 0)  \,,\quad 
 v_2 = (0,\, \lambda) \,,\quad 
 v_3 = (-\lambda,\, 0) \,,\quad 
 v_4 = (0,\, -\lambda) \, .  
\monend 
We have now three  particle distributions:   $f \in {\rm D2Q5}$, 
$g_x \in {\rm D2Q4} $ and $g_y \in {\rm D2Q4} $. The natural question  
is to find an intrinsic method 
 to determine the   equilibrium values
 $ \,  f_j^{\rm eq}  \,  $ for $\, 0 \leq j \leq 4 \, $ 
and $ \,  \big(  g_{xj}^{\rm eq} \,, $ $ \, g_{yj}^{\rm eq} \big) \,   $ 
for  $\, 1 \leq j \leq 4 $.  
As in the one-dimensional case, a  key point is to be able to explicit the 
dual entropy. In this two-dimensional case, 
the entropy variables $ \, \varphi \in \R^3 \,$ can be written as 
%
  \moneqstar    \label{3-variable-entropique} 
  \varphi \, = \, ( \theta \,, u \,, v ) \, , \,\,
  \theta = {{\partial \eta}\over{\partial \rho}} 
    =  \smash { {{c^2}\over{\gamma - 1}} - {{u^2+v^2}\over{2}} } \, . 
  \monendstar 
We have now as suggested in (\ref{3-st-venant-edp}): 
%
 \moneqstar 
 \eta^* (\theta , \, u,\,  v) \equiv p \equiv 
  {{\rho_0}\over{2\, c_0^2}} \, \Big( \theta + {1\over2}  \big( u^2 + v^2 \big) \Big)^2 \,. 
 \monendstar 
In order to determine the equilibrium distributions, we search 
 convex  functions $ \, h_j^* (\theta , \, u,\,  v) \,$ 
for  $\, 0 \leq j \leq 4 \,  $ such that the first set of Perthame-Bouchut 
conditions (\ref{1-perthame-bouchut-hyp}) are satisfied: 
\moneq   \label{3-perthame-bouchut-1} 
\sum_{j=0}^4 \,  h_j^* (\theta , \, u,\,  v) \equiv 
\eta^* (\theta , \, u,\,  v) \, . 
\monend 
Then
%
\moneqstar   
f_j^{\rm eq}  = {{\partial  h_j^*}\over{\partial \theta}} \,, \quad  
g_{xj}^{\rm eq} = {{\partial  h_j^*}\over{\partial u}} \,, \quad  
g_{yj}^{\rm eq} = {{\partial  h_j^*}\over{\partial v}} \, . 
\monendstar 

\monitem 
We have also to take into account the dual entropy fluxes
 $ \, \zeta_\alpha \, $ in order to represent correctly the first order terms
of the model (\ref{1-hyp-N}) or (\ref{3-st-venant-edp}) in our case. 
With   second set of Perthame-Bouchut 
conditions (\ref{1-perthame-bouchut-hyp}): 
\moneq   \label{3-perthame-bouchut-2} 
 \sum_{j=0}^4 \,  v_j^1 \, h_j^* (\theta , \, u,\,  v)  \equiv \eta^* \, u  \,, \quad   
\sum_{j=0}^4 \,  v_j^2 \, h_j^* (\theta , \, u,\,  v)  \equiv \eta^* \, v  \, . 
\monend 
%
%
%
For  the D2Q5 stencil, the  conditions of (\ref{3-perthame-bouchut-1})
(\ref{3-perthame-bouchut-2}) take the form
\moneq   \label{3-st-venant-contraintes-2} 
 h_0^* +  h_1^* + h_2^* + h_3^* + h_3^*  \equiv p   \,, \quad 
  \lambda \, \big( h_1^* -   h_3^* \big) \,  \equiv p \, u  \,, \quad 
    \lambda \, \big( h_2^*  -   h_4^*  \big) \equiv p \, v  \, . 
\monend 
We mimic for shallow water in two space dimensions  what we have done for 
the one-dimensional case (\ref{2-potentiels-st-venant-1d}) 
and we suggest here  to set
%
\moneqstar 
h_0^*(\theta) =  {{a}\over{2}} \, K \, \theta^2 
\monendstar 
as previously. 
%
%
Because this function $ \, h_0 ^* \,$ does not depend explicitly
on the variables $u$ and $v$, we are not defining 
a D1Q5Q5Q5 scheme but simply a  D1Q5Q4Q4 vectorial lattice Boltzmann scheme! 
The positive parameter $a$ has to be fixed at best.  
Nevertheless, we have still a lot of degrees of freedom. 
We suggest moreover to  cut into two part the first  relation 
of (\ref{3-st-venant-contraintes-2}): 
\moneq   \label{3-st-venant-enplus} 
  h_1^* +   h_3^* \,=\,  {1\over2} \big( p - h_0^* \big)  \,, \quad 
   h_2^* +   h_4^*  \,=\,  {1\over2}  \big( p - h_0^* \big) \, . 
\monend 
We have now a set of  5 independent  equations 
(\ref{2-potentiels-st-venant-1d}), 
(\ref{3-st-venant-contraintes-2}) and 
(\ref{3-st-venant-enplus}) with 5 unknowns $ \, h^*_j \,$. 
The end of the algebraic resolution of the system 
(\ref{2-potentiels-st-venant-1d}), 
(\ref{3-st-venant-contraintes-2}) and 
(\ref{3-st-venant-enplus}) is completely elementary.  

 \monitem
When the potentials  $ \, h^*_j \,$ are known, the 
computation of the equilibrium values is easy. 
%
With the 5+4+4=13 particle distributions, we can construct  13 moments
for the  D2Q5Q4Q4  lattice Boltzmann scheme. We suggest the following
5 moments associated with the distribution $\, f_j \,$:
%
\moneqstar \left\{ \begin{array} {l}  
 \displaystyle   \rho  \, = \,  f_0 + f_1 + f_2 + f_3 + f_4      \,, \quad 
 J_{x, \, \rho} \,=\, \lambda \, \big( f_1 - f_3 \big)    \,, \quad   
 J_{y, \, \rho} \,=\, \lambda \, \big( f_2 - f_4 \big)  \,, \quad  
 \\ \displaystyle  
 \varepsilon_{\rho} \,= \, f_1 + f_2 + f_3 + f_4 - 4 \, f_0  \,, \quad   
 XX_{\rho} \,=\, f_1 - f_2 + f_3 - f_4  \, .
\end{array} \right. \monendstar   
For the 8 moments relative to the distributions  $\, g_{xj}  \,$ 
and  $\, g_{yj}  \,$, we have chosen 
%
 \moneqstar \left\{ \begin{array} {ll}  \displaystyle  
   q_{x} \,=\,  g_{x1} +  g_{x2} +  g_{x3} +  g_{x4}  \,, \quad   
   & \displaystyle  f_{xx} \,=\, \lambda \, \big( g_{x1}  -   g_{x3}  \big)  \,,    
   \\ \displaystyle  
   f_{xy} \,=\, \lambda \, \big( g_{x2}  -   g_{x4}  \big)    \,, \quad   
   & \displaystyle  XX_{u} \,=\,  g_{x1} -  g_{x2} +  g_{x3} -  g_{x4}  
\end{array} \right. \monendstar     
and
%
 \moneqstar \left\{ \begin{array} {ll}  \displaystyle 
   q_{y} \,=\,  g_{y1} +  g_{y2} +  g_{y3} +  g_{y4}   \,, \quad     
   & \displaystyle f_{yx} \,=\, \lambda \, \big( g_{y1}  -   g_{y3}  \big)   \,, 
   \\ \displaystyle    
   f_{yy} \,=\, \lambda \, \big( g_{y2}  -   g_{y4}  \big)    \,, \quad    
   & \displaystyle XX_{v} \,=\,  g_{y1} -  g_{y2} +  g_{y3} -  g_{y4}   \, .  
 \end{array} \right. \monendstar     

 \monitem 
The value at equilibrium of the previous moments 
can be explicited, taking into account 
that the three moments $\, \rho, $ $\, q_x \,$ and $ \, q_y \, $ are at equilibrium.
We have:
%
  \moneqstar \left\{ \begin{array} {ll}  \displaystyle 
     J_{x, \, \rho}^{\rm eq} \,=\, \rho \, u   \,=\, q_x   \,, \quad   
     & \displaystyle  J_{y, \, \rho}^{\rm eq} \,=\, \rho \, v   \,=\, q_y  \,, 
    \\ \displaystyle    
    \varepsilon_{\rho}^{\rm eq} \,=\, 
    \Big( 1 - {{5\, a}\over{2}} \Big) \, \rho 
    +  {{5}\over{4}} \, {{\rho_0 \, (u^2 + v^2 )}\over{c_0^2}}  \,, \quad    
    & \displaystyle  XX_{\rho}^{\rm eq} \,=\, 0 \, . 
 \end{array} \right. \monendstar     
We have also 
%
  \moneqstar \left\{ \begin{array} {lll}  \displaystyle 
    f_{xx}^{\rm eq} \,=\, \rho \, u^2 + p \,, 
    & \displaystyle  f_{xy}^{\rm eq} \,=\, \rho \, u \, v  \,, 
    & \displaystyle XX_{u}^{\rm eq} \,=\, 0 
    \\ \displaystyle    
   f_{yx}^{\rm eq} \,=\, \rho \, u \, v \,,  
   & \displaystyle f_{yy} \,=\, \lambda \, \big( g_{y2}  -   g_{y4}  \big) \,,   
   & \displaystyle  XX_{v}^{\rm eq} \,=\, 0 \, . 
\end{array} \right. \monendstar     
The multiple relaxation time algorithm can be implemented without difficulty.
It is just necessary to write a relation of the type (\ref{2-relaxation}) 
for the 10 moments that are not at equilibrium. 
Our present choice is the BGK variant of the scheme, 
with all parameters  $ \, s_\ell \, $  supposed  to be equal.
The boundary conditions of wall constraint, supersonic inflow or
supersonic outflow are treated with an easy adaptation of  the usual methods 
of bounce-back and ``anti-bounce-back''.

\bigskip \bigskip  \noindent {\bf \large 4) \quad First test cases} 

%
\monitem 
We propose two bidimensional  test cases for the shallow water equations:
 a stationary shock reflexion 
and a classical unstationary forward facing step first proposed by Emery~\cite {Em68}
 for gas dynamics. 
The first test case is a the reflexion of an incident shock wave of 
angle $ \, -\pi / 4 \,$ issued from a ``left'' state
into a new shock of angle 
$ \, {\rm atan} (4/3) \,$  due to the physical nature of the 
``top'' state (in green on the left  picture of Fig.~2) 
and the ``right'' state (in indigo). 
The exact solution is determined through the use of Rankine Hugoniot relations. 
We have chosen 
%
  \moneqstar \left\{ \begin{array} {lll}  \displaystyle 
  \rho_\ell \,=\, 1 \,, 
  & \displaystyle  u_\ell   \,=\, 1.59497132403753 \,, 
  & \displaystyle  v_\ell    \,=\, 0 \,, 
  \\ \displaystyle    
  \rho_t \,=\, 1.17150636388320 \,, 
  & \displaystyle  u_t \,=\, 1.47822089880855 \,, 
  & \displaystyle v_t \,=\,  -0.116750425228984 \,, 
  \\ \displaystyle   
  \rho_r \,=\,   1.38196199044604 \,, 
  & \displaystyle u_r \,=\,  1.33228286727232 \,, 
  & \displaystyle v_r   \,=\, 0 \,. 
 \end{array} \right. \monendstar     
%
The stationary result of the vectorial lattice Boltzmann scheme for this first test case
can be compared with the pure finite volume approach with the Godunov \cite {Go59} scheme
solving a discontinuity at each interface at each time step. 
We have used three meshes of\br
 35 $\times$ 20,  70 $\times$ 40 and  140 $\times$ 80
grid points. The iso-values of density are presented on Fig.~2.
The numerical results are similar.

\smallskip \smallskip \newpage 
\centerline { 
{ \includegraphics[width=.32  \textwidth]
 {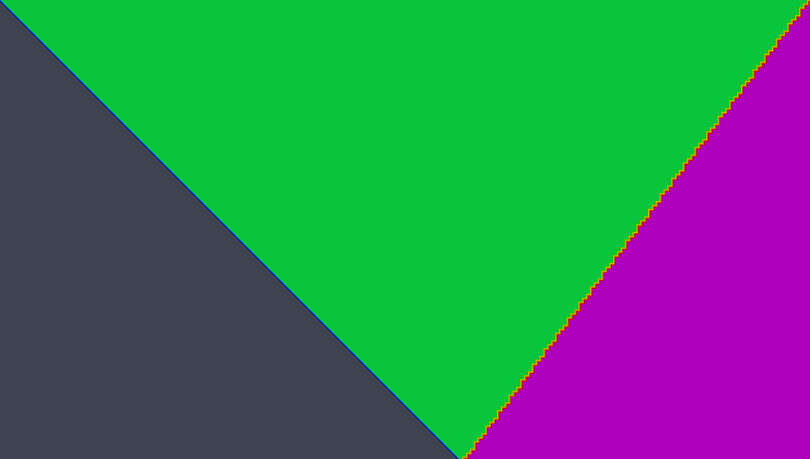}} 
{ \includegraphics[width=.32 \textwidth]
 {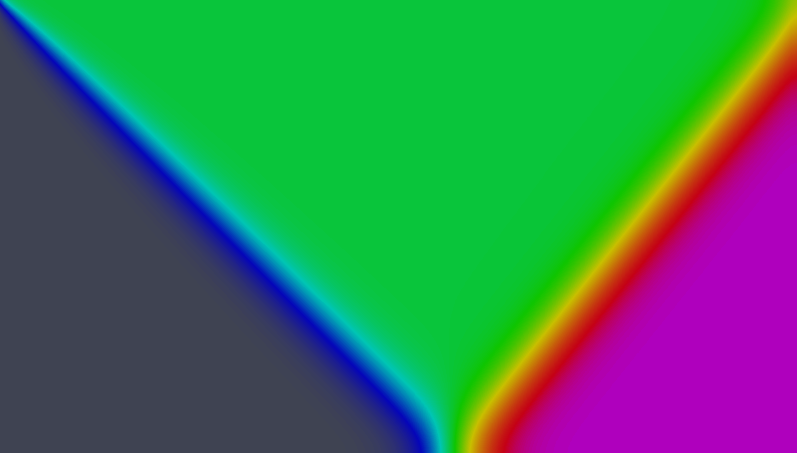}}   
{ \includegraphics[width=.32 \textwidth] 
{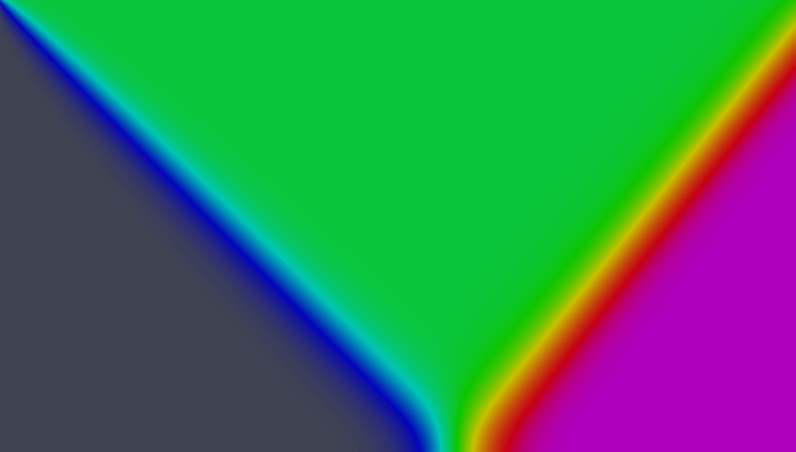}} } 

\noindent  {\bf Figure 2}. \quad 
Shock reflexion, mesh 140  $\times$ 80. 
Exact solution  (left), 
Lattice Boltzmann scheme D2Q5Q4Q4 (middle) and Godunov scheme (right). 
\smallskip  \smallskip 

\bigskip 

 \monitem
The second test case (Emery \cite {Em68}) is purely unstationary.
At time equal zero, a small step is created inside a  flow 
at Froude number equal to 3. A strong shock wave separates from the wall and various
nonlinear waves occur and interact. Our present experiment  (Fig.~3 and~4) 
shows  the ability of a vectorial lattice Boltzmann scheme to approach such a flow. 
We have refined the mesh, using three  families of meshes:
120 $\times$ 40, 240 $\times$ 80 and  480 $\times$ 120. 
We have used 
$ \, \lambda  =  80 $, 
$ \, a  =  0.05  $,
$ \, s_j  =  1.8 $ for all $j$ 
to  achieve  experimental stability. 
%
%
The time step is very small (due to the high value of $ \, \lambda = {{\Delta x}\over
{\Delta t}}$) and in consequence the computation relatively slow.

\bigskip  \bigskip  \centerline { 
{\includegraphics[width=.84 \textwidth] {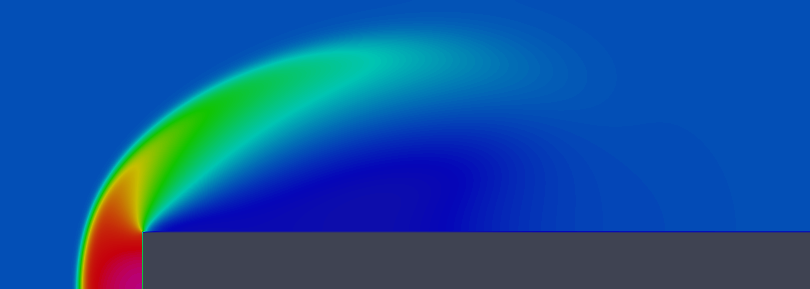}}}

\smallskip

\centerline { 
{\includegraphics[width=.84 \textwidth] {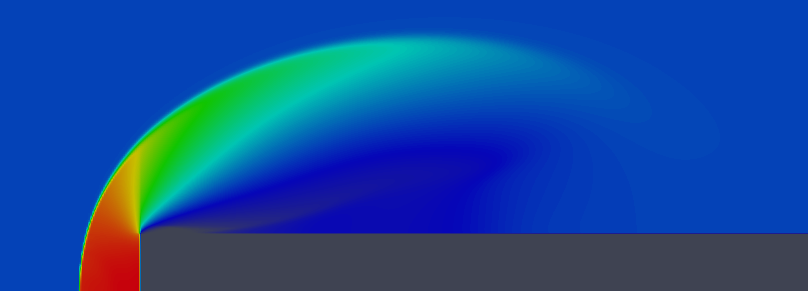}}}

\noindent  {\bf Figure 3}. \quad 
Emery test case for the shallow water equations, mesh  480 $\times$ 120,\br
 $ \, t  = 1/2 $,  density profile, 
D2Q5Q4Q4 vectorial lattice Boltzmann scheme (top) and  Godunov scheme (bottom). 
\smallskip  \smallskip 

\smallskip \smallskip \newpage 
\centerline { 
{\includegraphics[width=.84 \textwidth] {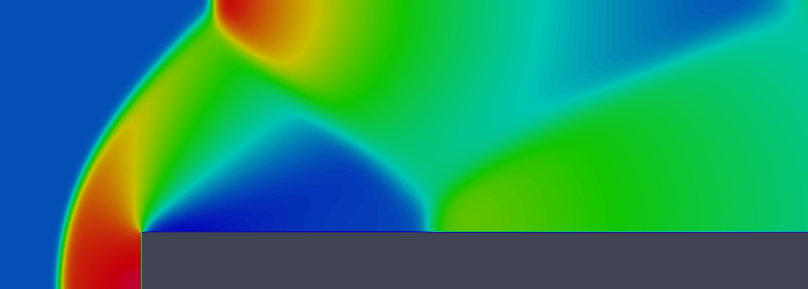}}}

\smallskip

\centerline { 
{\includegraphics[width=.84 \textwidth] {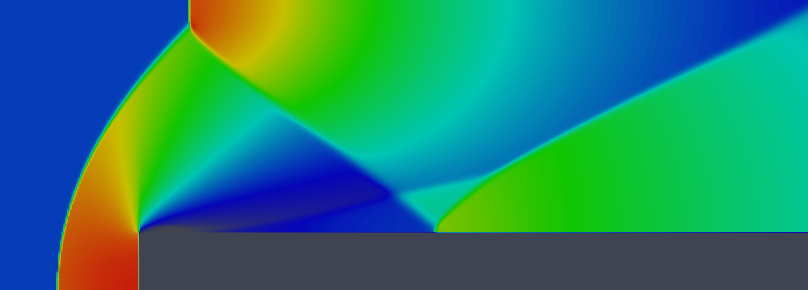}}}

\noindent  {\bf Figure 4}. \quad 
Emery test case for the shallow water equations, mesh  480 $\times$ 120,\br
 $ \, t = 4$,  density profile, 
D2Q5Q4Q4 vectorial lattice Boltzmann scheme (top) and  Godunov scheme (bottom). 
\smallskip   \smallskip 

\monitem  
We present our results for the finer mesh, at adimensionalized  time
equal to $1/2$ (Fig.~3) and $4$ (Fig.~4). 
The results show the ability for the vectorial scheme based on the decomposition 
of the dual entropy to capture such flows. 
Nevertheless, the Godunov scheme, well known of being only of order one, gives better 
unstationary result compared to the new approach.

\bigskip \bigskip  
\noindent {\bf \large  Conclusion}  

\monitem 
We have extended the methodology of  kinetic decomposition 
of the dual entropy previously studied for 1D problems  
into a general framework of vectorial lattice Boltzmann schemes 
for systems of conservation laws in several space dimensions, 
in the spirit of Bouchut \cite{Bo03}. 
The key point is to decompose the dual entropy of the system
into convex potentials satisfying the Perthame-Bouchut hypothesis.
Our first choices show that the system of shallow water equations
can be solved numerically without major difficulty.
Nevertheless, our first numerical  experiments  show that the
resulting scheme contains a lot of numerical viscosity.
Future work is  necessary to reduce this effect. 


\bigskip \bigskip  
\noindent {\bf \large  Acknowledgments}  

\monitem 
The author  thanks  Fran\c cois Bouchut for an enlightening
discussion during the elaboration of this  work.



\begin{thebibliography}{}

\end{thebibliography}


\bigskip \bigskip    
\noindent {\bf \large  References } 

 \vspace{-.4cm} 
 


\begin{thebibliography}{99}


\bibitem{ACCD92}
F. Alexander, H. Chen, S. Chen, G.  Doolen. 
  ``A lattice Boltzmann model for compressible fluids'', 
 {\it Phys.  Rev. A},  vol.~{\bf 46},   p~1967-1970, 1992.

\bibitem{BL87} 
B. Boghosian, C. Levermore. 
  ``A Cellular Automaton for Burgers's Equation'', 
  {\it  Complex Systems},  vol.~{\bf 1}, p~17-30, 1987.
 
\bibitem{BLY04} 
B. Boghosian, P. Love, J. Yepez.  
  ``Entropic Lattice Boltzmann Model for  Burgers' Equation'',   
{\it  Phil.  Trans. of the Royal Society A},  
vol.~{\bf 362},  p~1691-1702, 2004.
 

\bibitem {Bo99} 
F. Bouchut.  
 ``Construction of BGK models with a family of kinetic 
  entropies for a given system of conservation laws'', 
{\it Journal of  Statistical Physics}, vol.~{\bf 95}, p.~113-170, 1999.  

\bibitem {Bo03} 
F. Bouchut.  
  ``Entropy satisfying flux vector splittings and kinetic BGK models'', 
{\it Numerische  Mathematik},  vol.~{\bf  94},   p.~623-672, 2003.

\bibitem{CK09} 
S.  Chikatamarla, I. Karlin. 
  ``Lattices for the lattice Boltzmann method'', 
{\it Physical Review E}, vol.~{\bf 79}, 046701, 2009.

 
\bibitem{De02}  
P. Dellar. 
 ``Lattice Kinetic Schemes for Magnetohydrodynamics'', 
{\it Journal of Computational Physics},   vol.~{\bf 179},  p.~95-126, 2002.


\bibitem{DD05}  
 B. Despr\'es, F. Dubois. {\it  Syst\`emes hyperboliques de lois de conservation~;
 Application \`a  la dynamique des gaz}, 
  Editions de l'Ecole Polytechnique, Palaiseau, 2005. 


\bibitem{Du13}       
F. Dubois. 
  ``Stable lattice Boltzmann schemes with a dual entropy approach 
  for monodimensional nonlinear waves'',  
  {\it Computers and Mathematics with Applications}, 
 vol.~{\bf 65}, p.~142-159, 2013.



\bibitem{ELR93}    
B. Elton, C. Levermore, G.  Rodrigue. 
  ``Convergence of Convective-Diffusive Lattice Boltzmann Methods'', 
{\it SIAM  J. on Numerical  Analysis},  vol.~{\bf 32}, p.~1327-1354, 1995.



\bibitem{Em68}  
A.E. Emery.
  ``An evaluation of several differencing methods for inviscid fluid flow problems'', 
  {\it  Journal of Computational Physics},  vol.~{\bf 2}, p.~306–331, 1968. 


\bibitem{FL71}    
K.O. Friedrichs, P.D. Lax. 
  ``Systems of Conservation Equations with a Convex  Extension'', 
{\it Proc. Nat. Sciences USA}, vol.~{\bf 68}, p.~1686-1688, 1971. 



\bibitem{Go59}
S.K. Godunov.
  ``A Difference Scheme for Numerical Solution of Discontinuous 
  Solution of Hydrodynamic Equations'',  
{\it Math. Sbornik},  vol.~{\bf 47}, p.~271-306, 1959. 



\bibitem{DDH92}    
D. d'Humi\`eres. 
  ``Generalized Lattice-Boltzmann Equations'', 
in {\it Rarefied Gas Dynamics: Theory
and Simulations}, vol.~{\bf 159} of {\it AIAA Progress in
Astronautics and Astronautics}, p.~450-458, 1992.  


\bibitem{KA10}
I.V. Karlin, P. Asinari. 
 ``Factorization symmetry in the lattice Boltzmann method'',
{\it Physica A},   vol.~{\bf  389},  p.~1530-1548, 2010. 

\bibitem{KGSB98}
I.V.  Karlin, A.N. Gorban,  S. Succi and V. Boffi. 
  ``Maximum Entropy Principle for Lattice Kinetic Equations'',  
{\it  Physical Review Letters},  vol.~{\bf 81}, p.~6-9, 1998. 



\bibitem{KP94}
B.  Khobalatte, B. Perthame. 
  ``Maximum principle on the entropy and second-order kinetic schemes'', 
{\it Mathematics of Computation},   vol.~{\bf 62}, p.~ 119-131, 1994. 

\bibitem{LL00}  
P. Lallemand, L-S. Luo. 
  ``Theory of the lattice Boltzmann method: 
   Dispersion, dissipation, isotropy, Galilean invariance, and stability'',
{\it Physical Review E}, vol.~{\bf 61}, p.~6546-6562, June 2000.  
 

\bibitem{NSC08}  
X. Nie, X. Shan, H.  Chen. 
  ``Thermal lattice Boltzmann model for gases with internal degrees of freedom'', 
 {\it Physical  Review E},  vol.~{\bf  77}, p.~035701(R), 2008.


\bibitem{PHSS90}  
P.C. Philippi, L.A. Hegele, R. Surmas, D.N. Siebert. 
  ``From the Boltzmann to the Lattice-Boltzmann equation: beyond BGK collision models'', 
{\it International Journal of Modern Physics C},  vol.~{\bf 18},  p.~556-565, 2007.  


\bibitem{SC93}
X. Shan and H. Chen. 
 ``Lattice Boltzmann model for simulating flows  with multiple phases and components'', 
{\it Physical  Review E}, 
 vol.~{\bf 47}, p.~1815–1819, 1993.  



\bibitem{WWLL13}
J. Wang, D. Wang, P. Lallemand, L-S. Luo. 
 ``Lattice Boltzmann simulations of thermal convective flows in two dimensions'', 
  {\it Computers and Mathematics with Applications}, 
  vol.~{\bf 65}, p.~262-286, 2013. 


\end{thebibliography}
\end{document}